\documentclass[oneside,reqno,a4paper,12pt]{amsart}
\usepackage{preamble}

\title{Differential modules: a perspective on Bass' question}
\author{David Nkansah}
\date{\today}

\begin{document}

\begin{abstract}
    Guided by the \(Q\)-shaped derived category framework introduced by Holm and J{\o}rgensen, we provide a differential module analogue of a classical result that characterises when a finitely generated module over a local commutative noetherian ring has finite injective dimension. As an application, we characterise local Cohen–Macaulay rings using the homological algebra of differential modules.
\end{abstract}

\maketitle


{
\singlespacing
\hypersetup{linkcolor=black}
\tableofcontents
}

\section*{Introduction}
    
    A differential module is a module equipped with a square-zero endomorphism called a differential. The notion first appeared in the foundational book by Cartan and Eilenberg \cite{cartan-homological-1956}*{Chp.\ IV(1)}, where differential modules served primarily as a stepping stone towards the study of complexes. A complex can be viewed as a differential module endowed with a grading such that the differential has degree one. Accordingly, Cartan and Eilenberg defined cocycles, boundaries, and cohomology for differential modules before specialising these notions to graded settings.
    
   In recent years, differential modules have gained attention in commutative algebra and representation theory (see, for example, \cites{avramov-class-2007,rouquier-representation-2006,wei-gorenstein-2015,ringel-representations-2017,stai-differential-2017,brown-minimal-2022,holm-minimal-2025}).
    
    This paper is driven by the following three questions:
    \begin{itemize}
        \item [1.] Is there a suitable theory of \textit{homological algebra without grading}?
        \item [2.] How would such a theory compare to the homological algebra of complexes?
        \item [3.] Could this theory detect key properties of rings?
    \end{itemize}
    
   \paragraph{\bf For question one:} an affirmative answer is implicit in recent work by Holm and J{\o}rgen\-sen, who introduced the \(Q\)-shaped derived category \cite{holm-qshaped-2022} (see also \cite{holm-qshaped-2024}), a generalisation of the classical derived category. Their framework enables systematic homological constructions across various \(Q\)-shaped settings---one of which is the category of differential modules. We like to think of this as a theory of \textit{homological algebra without grading}.
   
   For example, (minimal) injective resolutions of modules appear naturally as complexes and motivate the notion of (minimal) semi-injective complexes. However, these resolutions do not directly suggest what a (minimal) semi-injective differential module should be. The \(Q\)-shaped framework developed in \cites{holm-qshaped-2022,holm-minimal-2025} addresses this gap and provides access to such notions.
   
   Moreover, differential modules admitting a projective flag---as considered in \cite{avramov-class-2007}*{2.1 and 2.2}---turn out to be a special case of the semi-projective differential modules arising naturally within the \(Q\)-shaped theory (see Example~\ref{ex: projective flag implies semi-projective}). We therefore situate the category of differential modules within the \(Q\)-shaped framework.
   
   \paragraph{\bf For question two:} we explore the adjunction triple \((\comp,\exp,\cocomp)\) (see \eqref{eqn: adjunction triple exp and comp}) between the category of differential modules and the category of complexes. The adjunction pair \((\comp,\exp)\) was previously considered in \cite{avramov-class-2007}*{pg.\ 7}. It behaves well on projective-minded objects---reflecting the authors' focus on this adjunction pair and differential modules admitting a projective flag.
   
   Recently in \cite{holm-minimal-2025}, a theory of minimal semi-injective \(Q\)-shaped objects was introduced. This motivates considering the remaining adjunction pair \((\exp,\cocomp{})\) into the picture, as it interacts nicely with injective-minded objects. We do so in this paper.

   \paragraph{\bf For question three:} homological algebra underwent a major breakthrough when it was used to characterise when a local commutative noetherian ring is regular. This was done by Serre \cite{serre-dimension-1956} and independently by Auslander and Buchsbaum \cite{auslander-homological-1956}. Later, Bass asked whether the existence of a finitely generated module of finite injective dimension suffices to characterise local Cohen–Macaulay rings (see \cite{bass-ubiquity-1963}*{Rmk.\ on pg.\ 14}). Using classical homological methods, this question was resolved---first in positive equicharacteristic by Peskine and Szpiro \cite{peskine-dimension-1973}, and then in full generality by Roberts \cites{roberts-theoreme-1987,roberts-multiplicities-1998}.
   
   Inspired by these works, we offer a differential modules analogue of this classical result. The key insight is as follows:

   Consider a finitely generated module \(M\) over an arbitrary ring. Viewing \(M\) as a stalk complex, we let \(I\) be its minimal semi-injective resolution. Then \(M\) has finite injective dimension if and only if \(I\) has finite support (i.e.\ \(I^i = 0\) for all but finitely many integers \(i\)).
   
   For differential modules, the situation is more subtle: since a differential module always has finite support, the standard finiteness condition does not suffice. However, in the setting of a local noetherian ring, we observe that a minimal semi-injective resolution has finite support if and only if the injective envelope of the residue field appears only finitely often. We prove that this characterisation remains valid in the category of differential modules.

\subsection*{The main theorem}

    Our main result gives a differential module analogue of a classical criterion for finitely generated modules of finite injective dimension. The result traces back to \cite{bass-ubiquity-1963}*{Lem.\ 2.7} (a complete modern reference is \cite{christensen-derived-2024}*{Thm.\ 16.4.8}). As a corollary, we provide a characterisation of local Cohen–Macaulay rings using differential modules.

   \begin{ThmIntro}[{Theorem~\ref{thm: finite injective dimension IFF bass number finite} and Corollary~\ref{thm: Bass' question}}]
       Let \((R,\fkm,k)\) be a local commutative noetherian ring and let \(M\) be finitely generated \(R\)-module. Consider the differential \(R\)-module \((M,0)\) with zero differential and choose a minimal semi-injective resolution \((M,0)\xrightarrow[]{}(I,d_I)\). Then the following statements are equivalent:
       \begin{itemize}
           \item[1.] \(M\) has finite injective dimension.
           \item[2.] The number of direct summands of \(I\) that are isomorphic to the injective envelope of the residue field \(k\) is finite.
       \end{itemize}
       Furthermore, \(R\) is Cohen-Macaulay if and only if there exists a nonzero finitely generated \(R\)-module \(M\) satisfying part 1 or part 2.
   \end{ThmIntro}

   The proof of this theorem hinges on the fact that the rightmost adjoint \(\cocomp\) in the triple adjunction \((\comp,\exp,\cocomp)\) (see \eqref{eqn: adjunction triple exp and comp}) preserves semi-injective complexes. This prompts a closer examination of the structural properties of this functor. In particular, we investigate when \(\cocomp\) preserves \textit{minimality}. We identify three instances of an algebra \(A\) together with a subclass of minimal semi-injective complexes
    \[
        \srC_A \subseteq \{\text{minimal semi-injective complexes of \(A\)-modules}\},
    \]
    for which minimality is preserved under \(\cocomp\). These cases are summarised in the following table:
    \begin{table*}[h]
    \centering
    \begin{tabularx}{\textwidth}{|>{\centering\arraybackslash}X|>{\centering\arraybackslash}c|>{\centering\arraybackslash}c|}
    \hline
    \textbf{Prop.} & \textbf{\(A\)} & \textbf{\(\srC_A\)} \\ \hline
    \ref{prop: cocompressing minimal semi-injectives of length 2} & Arbitrary & Complexes of length 2 \\ \hline
    \ref{prop: cocompressing minimal semi-injectives of fin gen} & Artinian & Bounded complexes \\ \hline
    \ref{prop: cocomp preserves local noetherian krull dimension one} & Local commutative noetherian of Krull dimension 0 & Arbitrary \\ \hline
    \end{tabularx}
    \end{table*}

\subsection*{Structure of the paper and the global setup}
    \paragraph{\bf Structure of the paper.} In Section~\ref{sec: 1. A lemma}, we prove a general lemma concerning the preservation of homologically defined objects under adjunctions arising in the \(Q\)-shaped framework. Section~\ref{sec: exp and (co)comp} develops an adjunction triple between the category of differential modules and the category of complexes. We then apply the previous lemma to study the behaviour of this adjunction triple. In Section~\ref{sec: diff tensor and hom}, we define an ad hoc tensor product and internal Hom in the category of differential modules and study their interaction with the adjunctions. Section~\ref{sec: bass numbers} specialises to the local noetherian case and uses the preceding theory to prove our main theorem. As a corollary, we obtain a differential module version of a classical result in commutative algebra. Section~\ref{sec: preservation of min semi-inj} establishes conditions under which the minimality of semi-injective objects is preserved under the functor \(\cocomp\) from Section~\ref{sec: exp and (co)comp}. Finally, in Section~\ref{sec: examples}, we provide examples illustrating the developed theory.

    \paragraph{\bf Global Setup.} We work in the following setup:
    \begin{itemize}
        \item \(\bbk\) is a hereditary commutative noetherian ring (for example,\ the ring of integers \(\bbZ\) or a field \(k\)). Hereafter, all tensor products with no subscripts are understood to be taken over this ground ring \(\bbk\).
        \item A \(\bbk\)-category is a category enriched over the closed symmetric monoidal category of \(\bbk\)-modules.
    \end{itemize}
    For a \(\bbk\)-algebra \(A\), we consider the following categories:
    \begin{itemize}
        \item \(\Mod{A}\) is the \(\bbk\)-category of right \(A\)-modules.
        \item \(\mod{A}\) is the full subcategory of \(\Mod{A}\) consisting of the finitely generated \(A\)-modules.
        \item \(\Ch{A}\) is the \(\bbk\)-category of (cohomologically graded) complexes of right \(A\)-modules.
        \item  \(\Dif{A}\) is the \(\bbk\)-category of differential right \(A\)-modules. Its objects are pairs \((X,d_X)\) consisting of an object \(X\) in \(\Mod{A}\) and an endomorphism \(X \xrightarrow[]{d_X} X\) in \(\Mod{A}\) such that the equation \(d_X d_X=0\) holds. We refer to such a morphism \(d_X\) a \textit{differential}. The Hom \(\bbk\)-module \(\Hom_{\Dif{A}}((X,d_X),({X'},d_{X'}))\) is the submodule of \(\Hom_{\Mod{A}}(X,{X'})\) consisting of those morphisms \(f\) such that the equation \(d_{X'}f = f d_X\) holds. In this case, we say \(f\) commutes or intertwines with the differentials. We sometimes represent a morphism \((X,d_X) \xrightarrow[]{f} (X',d_{X'})\) in \(\Dif{A}\) as a morphism between the corresponding differentials. That is, a commutative diagram
        \[
        \begin{tikzcd}
            X \arrow[r, "d_X"] \arrow[d, "f"'] & X \arrow[d, "f"] \\
            X' \arrow[r, "d_{X'}"]             & X'.              
        \end{tikzcd}
        \]
    \end{itemize}

    This paper uses the theory of semi-injective and semi-projective complexes. We suggest \cite{christensen-derived-2024}*{Sec.\ 5.2 and Sec.\ 5.3} as a reference.

\section{A lemma}\label{sec: 1. A lemma}

    The following lemma concerns abelian \(\bbk\)-categories of the following type: given a \(\bbk\)-algebra \(A\) and a small \(\bbk\)-category \(\clC\), we have the abelian \(\bbk\)-category 
    \[
    \Mod({\clC}\otimes A) \coloneqq \Fun_{\bbk}(\clC^{\op},\Mod{A})
    \]
    of \(\bbk\)-functors from \(\clC^{\op}\) to \(\Mod{A}\) (the notation adopted in \cite{holm-qshaped-2022} for the \(\bbk\)-category \(\Mod(\clC\otimes A)\) is \(\Mod_{\clC,A}\)). It is a locally finitely presentable Grothendieck abelian category (see \cite{holm-qshaped-2022}*{Prop.\ 3.12}) with all set indexed limits and set indexed colimits (see \cite{stenstrom-grundlehren-1975}*{Ch.\ V, preamble and Cor.\ X.4.4}): it has enough projective objects (see \cite{holm-qshaped-2022}*{Prop.\ 3.12(a)}), enough injective objects (see \cite{grothendieck-quelques-1957}*{Thm.\ 1.10.1}) and a finitely presented projective generator (see \cite{holm-qshaped-2022}*{Prop.\ 3.12(a)}). Its Hom \(\bbk\)-functor is denoted by \(\Hom_{\clC\otimes A}\) and the right derived functor of \(\Hom_{\clC\otimes A}\) is denoted by \(\Ext^1_{\clC\otimes A}\).
    When \(\clC\) satisfies the extra conditions given in \cite{holm-qshaped-2024}*{Setup~2.9}, the abelian \(\bbk\)-category \(\Mod(\clC\otimes A)\) permits the following notions:
    \begin{itemize}
        \item The notion of an `exact' object (see \cite{holm-qshaped-2022}*{Def.\ 4.1 and Thm.\ 7.1}). This generalises the acyclic (or exact) complexes.
    \end{itemize}
    The exact objects give rise to the following two classes of objects (which were also considered in \cite{holm-qshaped-2022}):
    \begin{itemize}
        \item An object \(P\) in \(\Mod(\clC \otimes A)\) is \textit{semi-projective} if \(\Ext^1_{\clC\otimes A}(P,E)=0\) for each exact object \(E\) in \(\Mod(\clC\otimes A)\).
        \item An object \(I\) in \(\Mod(\clC \otimes A)\) is \textit{semi-injective} if \(\Ext^1_{\clC\otimes A}(E,I)=0\) for each exact object \(E\) in \(\Mod(\clC\otimes A)\).
    \end{itemize}
    To parse the following lemma, the reader is invited to replace the categories involved with the category of complexes where the exact objects are the acyclic complexes.
    
\begin{lem}\label{lem: (F,R) adjunction preservation}
    Let \(A\) and \(B\) be \(\bbk\)-algebras and \(Q\) and \(\Gamma\) be small \(\bbk\)-categories that satisfy  \cite{holm-qshaped-2024}*{Setup~2.9}. Consider an adjunction pair \((L,R)\) as follows:
    \[
    \begin{tikzcd}
        \Mod(Q \otimes A) \arrow[r, "L", shift left] & {\Mod(\Gamma \otimes B),} \arrow[l, "R", shift left]
    \end{tikzcd}
    \]
    and assume \(L\) and \(R\) are exact functors. The following statements hold:
    \begin{itemize}
        \item[1.] If \(I\) is an injective object in \(\Mod(\Gamma \otimes B)\), then \(RI\) is an injective object in \(\Mod(Q \otimes A)\).
        \item[2.] If \(P\) is a projective object in \(\Mod(Q \otimes A)\), then \(LP\) is a projective object in \(\Mod(\Gamma \otimes B)\).
        \item[3.] Let \(X\) be an object in \(\Mod(Q \otimes A)\) and let \(Y\) be an object in \(\Mod(\Gamma \otimes B)\). There is a natural isomorphism of \(\bbk\)-modules
                \[
                \Ext^1_{\Gamma \otimes B}(LX,Y) \cong \Ext^1_{Q \otimes A}(X,RY).
                \]
        \item[4.] Suppose \(LE\) is an exact object in \(\Mod(\Gamma \otimes B)\) whenever \(E\) is an exact object in \(\Mod(Q \otimes A)\). Then \(RI\) is a semi-injective object in \(\Mod(Q \otimes A)\) whenever \(I\) is a semi-injective object in \(\Mod(\Gamma \otimes B)\).
        \item[5.] Suppose \(RE\) is an exact object in \(\Mod(Q \otimes A)\) whenever \(E\) is an exact object in \(\Mod(\Gamma \otimes B)\). Then \(LP\) is a semi-projective object in \(\Mod(\Gamma \otimes B)\) whenever \(P\) is a semi-projective object in \(\Mod(Q \otimes A)\).
    \end{itemize}
\begin{proof}
    \textit{Part 1.} The functor \(\Hom_{\Gamma \otimes B}(L(-),I)\) is exact as, by the assumptions, it is the composition of exact functors. The result follows by the adjunction isomorphism \(\Hom_{\Gamma \otimes B}(L(-),I) \cong \Hom_{Q \otimes A}(-,RI)\).

    \textit{Part 2.} Follows dually to part 1.
    
    \textit{Part 3.} Let \(I^{0} \xrightarrow{} I^{1} \xrightarrow{} I^{2}\) be a truncated injective resolution of \(Y\). By part~1 and as \(R\) is an exact functor, the sequence \(R(I^{0}) \xrightarrow{} R(I^{1}) \xrightarrow{} R(I^{2})\) is a truncated injective resolution of \(RY\). We apply the functor \(\Hom_{Q \otimes A}(X,-)\) to this shorted injective resolution and use the adjunction isomorphism of the adjunction pair \((L,R)\) to get the following commutative diagram of \(\bbk\)-modules
    \[
    \begin{tikzcd}
        {\Hom_{\Gamma \otimes B}(LX,I^0)} \arrow[r] \arrow[d, "\cong"] & {\Hom_{\Gamma \otimes B}(LX,I^1)} \arrow[r] \arrow[d, "\cong"] & {\Hom_{\Gamma \otimes B}(LX,I^2)} \arrow[d, "\cong"] \\
        {\Hom_{Q \otimes A}(X,R(I^0))} \arrow[r]                    & {\Hom_{Q \otimes A}(X,R(I^1))} \arrow[r]                    & {\Hom_{Q \otimes A}(X,R(I^2))},    
    \end{tikzcd}
    \]
    where the vertical arrows are the adjunction isomorphisms. As any choice of injective resolution computes \(\Ext\), the claim follows by taking cohomology.

    \textit{Part 4.} Let \(I\) be a semi-injective object in \(\Mod(\Gamma \otimes B)\). For each exact object \(E\) in \(\Mod(Q \otimes A)\), the object \(LE\) is exact by assumption and therefore the \(\bbk\)-module \(\Ext^1_{\Gamma \otimes B}(LE,I)\) is zero. The result follows by the isomorphism found in part~3.

    \textit{Part 5.} Follows dually to part 4.
\end{proof}
\end{lem}

\section{The expansion and (co)compression functors}\label{sec: exp and (co)comp}
    
    Consider a \(\bbk\)-algebra \(A\). Let \(\Gamma\) be the \(\bbk\)-category of paths of the quiver
\begin{equation}\label{quiver cpx}
\begin{tikzcd}
        \cdots \arrow[r] & 1 \arrow[r] & 0 \arrow[r] & -1 \arrow[r] & \cdots
\end{tikzcd}
\end{equation}
    modulo the relation that the concatenation of any two consecutive arrows is zero. The categories \(\Mod(\Gamma \otimes A)\) and \(\Ch{A}\) are isomorphic and we will henceforth identify them. Similarly, let \(Q\) be the \(\bbk\)-category of paths of the quiver
\begin{equation}\label{quiver dif}
\begin{tikzcd}
        * \arrow[loop, distance=2em, in=55, out=125]
\end{tikzcd}
\end{equation}
    modulo the relation that the concatenation of the unique loop with itself is zero. The categories \(\Mod(Q \otimes A)\) and \(\Dif{A}\) are isomorphic and we will henceforth identify them. Let \(\Gamma\xrightarrow[]{\pi}Q\) be the \(\bbk\)-functor induced by the unique morphism of quivers from the quiver \eqref{quiver cpx} to the quiver \eqref{quiver dif}. By the theory of Kan extensions (see \cite{maclane-categories-1998}*{Sec.\ X.3, Cor.\ 2 and its dual}), there is an adjunction triple \((\comp,\exp,\cocomp)\) as follows:
    \begin{equation}\label{eqn: adjunction triple exp and comp}
    \begin{tikzcd}
        \Dif{A} \arrow[rr, "\exp" description] &  & {\Ch{A}}, \arrow[ll, "\cocomp", bend left] \arrow[ll, "\comp"', bend right]
    \end{tikzcd}
    \end{equation}
    where \(\exp\) is given by precomposition by \(\pi\). Inspired by the terminology in \cite{avramov-class-2007}, where the top two functors appear, we call \(\exp\) the \textit{expansion functor}, \(\comp\) \textit{the compression functor} and \(\cocomp\) \textit{the cocompression functor}. For convenience, the explicit actions of these functors are given below.
    
    \paragraph{\bf Expansion functor.} For an object \((X,d_X)\) in \(\Dif{A}\), the expansion \(\exp(X,d_X)\) is the object in \(\Ch{A}\) given by the diagram
    \[
    \cdots \xrightarrow[]{} X \xrightarrow[]{d_X} X \xrightarrow[]{d_X} X \xrightarrow[]{} \dots
    \]
    and for a morphism \((X,d_X) \xrightarrow[]{f} (X',d_{X'})\) in \(\Dif{A}\), the expansion
    \[
    \exp(X,d_X) \xrightarrow[]{\exp{f}} \exp(X',d_{X'})
    \]
    is the morphism in \(\Ch{A}\) given by the commutative diagram
    \[
    \begin{tikzcd}
        \cdots \arrow[r] & X \arrow[r, "d_X"] \arrow[d, "f"] & X \arrow[r, "d_X"] \arrow[d, "f"] & X \arrow[r] \arrow[d, "f"] & \cdots \\
        \cdots \arrow[r] & X' \arrow[r, "d_{X'}"]                & X' \arrow[r, "d_{X'}"]                & X' \arrow[r]                & \cdots.
\end{tikzcd}
    \]

\paragraph{\bf Compression functor.} For an object
    \[
    Y = \cdots \xrightarrow[]{} Y^{i-1} \xrightarrow[]{d^{i-1}_Y} Y^{i}  \xrightarrow[]{d^{i}_Y} Y^{i+1} \xrightarrow[]{} \dots
    \]
    in \(\Ch{A}\), the compression \(\comp{Y}\) is the object in \(\Dif{A}\) given by the pair
    \[
    \left( \coprod_{i\in\bbZ} Y^i,\: \coprod_{i\in\bbZ} d^i_Y\right),
    \]
    where \(\coprod_{i\in\bbZ} d^i_Y\) is the canonical morphism induced by the compositions \(Y^j \xrightarrow[]{d^j_Y} Y^{j+1} \xrightarrow[]{} \coprod_{i\in \bbZ} Y^i\), for each integer \(j\), where \(Y^{j+1} \xrightarrow[]{} \coprod_{i\in \bbZ} Y^i\) are the canonical inclusions. For a morphism \(Y \xrightarrow[]{f} Z\):
    \[
    \begin{tikzcd}
        \cdots \arrow[r] & Y^{i-1} \arrow[r, "d^{i-1}_Y"] \arrow[d, "f^{i-1}"] & Y^{i} \arrow[r, "d^{i}_Y"] \arrow[d, "f^i"] & Y^{i+1} \arrow[r] \arrow[d, "f^{i+1}"] & \cdots \\
        \cdots \arrow[r] & Z^{i-1} \arrow[r, "d^{i-1}_Z"]                      & Z^{i} \arrow[r, "d^{i}_Z"]                  & Z^{i+1} \arrow[r]                      & \cdots
    \end{tikzcd}
    \]
    in \(\Ch{A}\), the compression \(\comp{Y} \xrightarrow[]{\comp{f}} \comp{Z}\) is the canonical morphism \(\coprod_{i\in\bbZ}f^i\) in \(\Dif{A}\), induced by the compositions \(Y^j \xrightarrow[]{f_j} Z^j \xrightarrow[]{} \coprod_{i\in\bbZ} Z^i\), where \(Z^j \xrightarrow[]{} \coprod_{i\in\bbZ} Z^i\) are the canonical inclusion.

\paragraph{\bf Cocompression functor.} For an object
    \[
    Y = \cdots \xrightarrow[]{} Y^{i-1} \xrightarrow[]{d^{i-1}_Y} Y^{i}  \xrightarrow[]{d^{i}_Y} Y^{i+1} \xrightarrow[]{} \dots
    \]
    in \(\Ch{A}\), the cocompression \(\cocomp{Y}\) is the object in \(\Dif{A}\) given by the pair
    \[
    \left( \prod_{i\in\bbZ} Y^i,\: \prod_{i\in\bbZ} d^i_Y\right),
    \]
    where \(\prod_{i\in\bbZ} d^i_Y\) is the canonical morphism induced by the compositions \(\prod_{i\in\bbZ} Y^i \xrightarrow[]{} Y^j \xrightarrow[]{d^j_Y} Y^{j+1}\), for each integer \(j\), where \(\prod_{i\in\bbZ} Y^i \xrightarrow[]{} Y^j\) are the canonical projections. For a morphism \(Y \xrightarrow[]{f} Z\):
    \[
    \begin{tikzcd}
        \cdots \arrow[r] & Y^{i-1} \arrow[r, "d^{i-1}_Y"] \arrow[d, "f^{i-1}"] & Y^{i} \arrow[r, "d^{i}_Y"] \arrow[d, "f^i"] & Y^{i+1} \arrow[r] \arrow[d, "f^{i+1}"] & \cdots \\
        \cdots \arrow[r] & Z^{i-1} \arrow[r, "d^{i-1}_Z"]                      & Z^{i} \arrow[r, "d^{i}_Z"]                  & Z^{i+1} \arrow[r]                      & \cdots
    \end{tikzcd}
    \]
    in \(\Ch{A}\), the cocompression \(\cocomp{Y} \xrightarrow[]{\cocomp{f}} \cocomp{Z}\) is the canonical morphism \(\prod_{i\in\bbZ}f^i\) in \(\Dif{A}\), induced by the morphism \(\prod_{i\in\bbZ} Y^i \xrightarrow[]{} Y^i \xrightarrow[]{f_i} Z^i\), where \(\prod_{i\in\bbZ}Y^i \xrightarrow[]{} Y^j\) are the canonical projections.
    
    We now recall some additional functors that will be used later (see \cite{holm-qshaped-2022}*{Cor.\ 3.9 and Prop.\ 7.15}).
    
    \paragraph{\bf Stalk functors.} These are functors 
    \[
        \Dif{A} \xleftarrow[]{\S_*} \Mod{A} \xrightarrow[]{\S_0} \Ch{A},
    \]
    where \(\S_*{M}=(M,0)\) and where \(\S_0\) identifies an \(A\)-module with its stalk complex in degree 0.
    
    \paragraph{\bf Evaluation functors and their left adjoint.} There are functors
    \[
        \Dif{A} \xrightarrow[]{\E_*} \Mod{A} \xleftarrow[]{\E_i} \Ch{A},
    \]
    where the index \(i\) ranges over the integers. The evaluations have object assignments given by \(\E_*(X,d_X)=X\) and \(\E_i(X^\bullet)=X^i\). The evaluations admit left adjoints
    \[
    \Dif{A} \xleftarrow[]{\F_*} \Mod{A} \xrightarrow[]{\F_i} \Ch{A}, 
    \]
    with object assignments given by
    \[
    \F_*{M} = \left(M\oplus M,\: {\begin{pmatrix}0 \amsamp 0 \\ \id_M \amsamp0 \end{pmatrix}}\right) \text{\qquad and \qquad}\F_i{M} =\cdots \xrightarrow[]{} 0 \xrightarrow[]{} M \xrightarrow[]{\id_M} M \xrightarrow[]{} 0 \xrightarrow[]{} \cdots,
    \]
    where \(F_i{M}\) is concentrated in cohomological degrees \(i\) and \(i+1\).

    \paragraph{\bf Cohomology functors.} These are functors
    \[
    \Dif{A} \xrightarrow[]{\Ho^*} \Mod{A} \xleftarrow[]{\Ho^i} \Ch{A},
    \]
    where the index \(i\) ranges over the integers. The cohomology functor \(\Ho^*\) has object assignments given by \(\Ho^*(X,d_X) = \Ker{d_X}/\Img{d_X}\) and \(\Ho^i\) is the usual \(i\)-th cohomology of an object in \(\Ch{A}\).

    We now address some definitions in the category \(\Dif{A}\) (compare these to the analogous definitions in \(\Ch{A}\)).

    \paragraph{\bf Acyclic differential modules and quasi-isomorphisms.} An object \((X,d_X)\) in \(\Dif{A}\) is \textit{acyclic} if \(\Ho^*(X,d_X)=0\). A morphism \((X,d_X) \xrightarrow[]{f} (Y,d_Y)\) in \(\Dif{A}\) is a \textit{quasi-isomorphism} if the obvious induced morphism \(\Ho^*(f)\) is an isomorphism in \(\Mod{A}\).

    \paragraph{\bf Contractible differential modules.} An object \((X,d_X)\) in \(\Dif{A}\) is \textit{contractible} if there exists an endomorphism \(X \xrightarrow[]{h} X\) in \(\Mod{A}\) such that the equation \(\id_X=hd_X +d_X h\) holds. We call \(h\) a \textit{contracting homotopy} of \((X,d_X)\).
    
    \begin{rem}\label{rem: all defn of contractible are the same}
    In \cite{ringel-representations-2017}, the notion of a morphism in \(\Dif{A}\) being homotopic to zero was considered. In their terminology, an object is contractible if its identity morphism is homotopic to zero. In \cite{avramov-class-2007}*{1.5 and Prop.\ 1.8}, alternative definitions of a contractible object are considered. Proposition~\ref{prop: contractible equivalence} shows that all these definitions are equivalent.
    \end{rem}
    
    \paragraph{\bf Minimal semi-injective resolutions.} An object \((I,d_I)\) in \(\Dif{A}\) is \textit{semi-injective} if for each acyclic object \((E,d_E)\) in \(\Dif{A}\), we have \(\Ext^1_{\Dif{A}}((E,d_E),(I,d_I))=0\) (this is a special case of the semi-injective objects defined in Section~\ref{sec: 1. A lemma}). A \textit{minimal semi-injective object} in \(\Dif{A}\) is a semi-injective object whose only injective subobject is the zero object. A \textit{semi-injective resolution} of an object \((X,d_X)\) in \(\Dif{A}\) is a quasi-isomorphism \((X,d_X) \xrightarrow[]{} (I,d_I)\) with \((I,d_I)\) a semi-injective object in \(\Dif{A}\). If, in this case, \((I,d_I)\) is a minimal semi-injective object, then we say \((X,d_X) \xrightarrow[]{} (I,d_I)\) is a \textit{minimal semi-injective resolution} of \((X,d_X)\) in \(\Dif{A}\).

    \begin{rem}
        The above definitions are just special cases of the definitions found in \cite{holm-minimal-2025}*{Def.\ A} applied to the category \(\Dif{A}\) and hence by \cite{holm-minimal-2025}*{Thm.\ 3.1}, we have access to equivalent characterisations of a minimal semi-injective object in \(\Dif{A}\).
    \end{rem}

    \paragraph{\bf Semi-projective objects.} An object \((P,d_P)\) in \(\Dif{A}\) is \textit{semi-projective} if for each acyclic object \((E,d_E)\) in \(\Dif{A}\), we have \(\Ext^1_{\Dif{A}}((P,d_P),(E,d_E))=0\) (this is a special case of the semi-projective objects defined in Section~\ref{sec: 1. A lemma}).

    \begin{ex}[Projective flags]\label{ex: projective flag implies semi-projective}
    In this example, we will suppress the differential of objects in \(\Dif{A}\). Let \(X\) be an object in \(\Dif{A}\). We say \(X\) \textit{admits a projective flag} (see \cite{avramov-class-2007}*{2.1 and 2.2}) if there is a nested sequence
    \[
    X^0 \subseteq X^1 \subseteq \cdots \subseteq X^{n} \subseteq \cdots
    \]
    of subobjects of \(X\) in \(\Dif{A}\) with the property that \(X^0=0\), \(X= \bigcup_{i\geq0}X^i\) and  for every \(i\geq0\), the quotients \(X^{i+1}/X^{i}\) are isomorphic to \(\S_*{P}\), for \(P\) a projective object in \(\Mod{A}\). Every object in \(\Dif{A}\) that admits a projective flag lies in the class of objects in \(\Dif{A}\) that have a filtration by objects in \(\{\S_*{P} \mid P \text{ is a projective object in }\Mod{A}\}\) (see \cite{holm-cotorsion-2019}*{6.4}). Therefore, by \cite{holm-qshaped-2024}*{Prop.\ 3.2(a)}, objects in \(\Dif{A}\) that admit a projective flag are also semi-projective. More generally, by the additivity of \(\Ext^1_{\Dif{A}}\), if \(X \xrightarrow[]{} Y\) is a retraction in \(\Dif{A}\) with \(X\) admitting a projective flag, then \(Y\) is also semi-projective.
\end{ex}

\begin{rem}
    Considering Example~\ref{ex: projective flag implies semi-projective} and noting that projective objects in \(\Dif{A}\) are in particular contractible (see \cite{cartan-homological-1956}*{Prop.\ 2.5 and Prop.\ 2.3}), the result \cite{avramov-class-2007}*{Thm.\ 2.3} is a special case of \cite{holm-qshaped-2022}*{Thm.\ 4.4(a)}. 
\end{rem}

    We prove some basic properties of the expansion and (co)compression functors.
    
\begin{prop}\label{prop: (co)comp are exact}
    Consider the adjunction triple \((\comp,\exp,\cocomp)\) as follows:
    \[
    \begin{tikzcd}
        \Dif{A} \arrow[rr, "\exp" description] &  & {\Ch{A}}. \arrow[ll, "\cocomp", bend left] \arrow[ll, "\comp"', bend right]
    \end{tikzcd}
    \]
    The following statements hold:
    \begin{itemize}
        \item[1.] The expansion functor \(\Dif{A} \xrightarrow[]{\exp} \Ch{A}\) is exact.
        \item[2.] The compression functor \(\Ch{A} \xrightarrow[]{\comp} \Dif{A}\) is exact.
        \item[3.] The cocompression functor \(\Ch{A} \xrightarrow[]{\cocomp} \Dif{A}\) is exact.
    \end{itemize}
\begin{proof}
    \textit{Part 1.} This follows as \(\exp\) is both a left and a right adjoint.
        
    \textit{Part 2.} As \(\comp\) is a left adjoint, it suffices to show that it is left exact. Let \(X \xrightarrow[]{f} Y\) be a monomorphism in \(\Ch{A}\). Equivalently, this means that for each integer \(i\), the morphism \(\E_if\) in \(\Mod{A}\) is a monomorphism. We have to show that \(\comp{X} \xrightarrow[]{\comp{f}}\comp{Y}\) is a monomorphism in \(\Dif{A}\) which occurs if and only if \(\E_*{\comp{X}} \xrightarrow[]{\E_*\comp{f}} \E_*\comp{Y}\) is a monomorphism in \(\Mod{A}\). But, \(\E_*\comp{f}\) is precisely \(\coprod_{i\in\bbZ}\E_i(f)\) which is a monomorphism as coproducts are exact in the AB4 category \(\Mod{A}\).

    \textit{Part 3.} This follows by a similar argument as seen in part 1, noting that products are exact in the AB4* category \(\Mod{A}\).
\end{proof}
\end{prop}

    Hence, the adjunction pairs \((\comp, \exp)\) and \((\exp, \cocomp)\) fit into Lemma~\ref{lem: (F,R) adjunction preservation}. We will now see that all working definitions of a contractible object in the literature are equivalent (see Remark~\ref{rem: all defn of contractible are the same}).

\begin{prop}\label{prop: contractible equivalence}
    Let \((X,d_X)\) be an object in \(\Dif{A}\). The following statements are equivalent:
    
    \begin{itemize}
        \item[1.] \((X,d_X)\) is a contractible object in \(\Dif{A}\).
        \item[2.] \(\exp(X,d_X)\) is a contractible object in \(\Ch{A}\).
        \item[3.] \((X,d_X)\) is isomorphic to \(\F_*{M}\) for some object \(M\) in \(\Mod{A}\).
        \item[4.] \((X,d_X)\) is acyclic and the short exact sequence
                \[
                   0 \xrightarrow[]{} \Ker{d_X} \xrightarrow{} X \xrightarrow[]{} \Img{d_X} \xrightarrow[]{} 0
                \]
                splits in \(\Mod{A}\).
    \end{itemize}
\begin{proof}
    The proof of statements 2, 3 and 4 being equivalent can be found in \cite{avramov-class-2007}*{Prop.\ 1.8}. We now prove the equivalence of statements 1 and 2.

    \((1 \Longrightarrow 2):\) Let \(X\xrightarrow[]{h} X\) be a contracting homotopy of \((X,d_X)\). Setting \(h^i=h\) for each integer \(i\), the sequence \((X \xrightarrow[]{h^i} X)_{i\in\bbZ}\) is a contracting homotopy for \(\exp(X,d_X)\).

    \((2 \Longrightarrow 1):\) Let \((X \xrightarrow[]{h^i} X)_{i\in\bbZ}\) be a contracting homotopy of \(\exp(X,d_X)\). Choose an integer \(i\) and set \(h\coloneqq h^{i+1}d_Xh^i\). The following equation \(\id_X = h^{i+1}d_X+d_Xh^i\) holds and implies both the equations \(d_X=d_Xh^{i}d_X\) and \(d_X=d_Xh^{i+1}d_X\). We have
    \begin{align*}
    \id_X = h^{i+1}d_X+d_Xh^i &= h^{i+1}(d_Xh^id_X) + (d_Xh^{i+1}d_X)h^i\\
                              &= (h^{i+1}d_Xh^i)d_X + d_X(h^{i+1}d_Xh^i) = hd_X+d_Xh. \qedhere
    \end{align*}
\end{proof}
\end{prop}

    The next result tells us how the cohomology functors interact with the expansion and (co)compression functors.

\begin{prop}\label{prop: cohomology of expansion and (co)homology}
    Let \((X,d_X)\) be an object in \(\Dif{A}\) and \(Y\) be an object in \(\Ch{A}\). The following statements hold:
    \begin{itemize}
        \item[1.] For each integer \(i\), we have an equality \(\Ho^i(\exp(X,d_X))=\Ho^*(X,d_X)\).
        \item[2.] There is a natural isomorphism \(\Ho^*(\comp{Y}) \cong \coprod_{i\in\bbZ}\Ho^i(Y)\) in \(\Mod{A}\).
        \item[3.] There is a natural isomorphism \(\Ho^*(\cocomp{Y}) \cong  \prod_{i\in\bbZ}\Ho^i(Y)\) in \(\Mod{A}\).
    \end{itemize}
    \begin{proof}
        \textit{Part 1.} This follows from the explicit description of the expansion functor \(\exp\).
        
        \textit{Part 2.} This is stated in \cite{avramov-class-2007}*{pg.\ 7} and a proof is similar to part 3, which we give below.
        
        \textit{Part 3.} We have a composition of natural isomorphisms
        \begin{align*}
        \Ho^*(\cocomp{Y}) = \Ker\left(\prod_{i\in\bbZ}d^i_Y\right) \Big/ \Img\left(\prod_{i\in\bbZ}d^i_Y\right) &\cong \prod_{i\in\bbZ}\Ker{d^i_Y} \Big/ \prod_{i\in\bbZ}\Img{d^i_Y}\\
                                                                                                                &\cong \prod_{i\in\bbZ}(\Ker{d^i_Y} / \Img{d^i_Y}) = \prod_{i\in\bbZ}\Ho^{i}(Y),
        \end{align*}
        that exist since products are exact in the AB4* category \(\Mod{A}\).
    \end{proof}
\end{prop}

\begin{cor}\label{cor: preserves acyclicity and quasi-iso}
    The following statements hold:
    \begin{itemize}
        \item[1.] \((X,d_X)\) is an acyclic object in \(\Dif{A}\) if and only if \(\exp(X,d_X)\) is an acyclic object in \(\Ch{A}\).
        \item[2.] \(Y\) is an acyclic object in \(\Ch{A}\) if and only if \(\comp{Y}\) is an acyclic object in \(\Dif{A}\).
        \item[3.] \(Y\) is an acyclic object in \(\Ch{A}\) if and only if \(\cocomp{Y}\) is an acyclic object in \(\Dif{A}\).
        \item[4.] \(f\) is a quasi-isomorphism in \(\Dif{A}\) if and only if \(\exp{f}\) is a quasi-isomorphism in \(\Ch{A}\).
        \item[5.] \(g\) is a quasi-isomorphism in \(\Ch{A}\) if and only if \(\comp{g}\) is a quasi-isomorphism in \(\Dif{A}\).
        \item[6.] \(g\) is a quasi-isomorphism in \(\Ch{A}\) if and only if \(\cocomp{g}\) is a quasi-isomorphism in \(\Dif{A}\). 
    \end{itemize}
\begin{proof}
    The statements follow from Proposition~\ref{prop: cohomology of expansion and (co)homology}.
\end{proof}
\end{cor}

    Our next aim is to study to which extent do the expansion and (co)compression functors interact with semi-injective and semi-projective objects.

\begin{prop}\label{prop: preserves injective}
    Let \((I,d_I)\) be an object in \(\Dif{A}\) and \(J\) be an object in \(\Ch{A}\). The following statements hold:
    \begin{itemize}
        \item[1.] If \(J\) is an injective object in \(\Ch{A}\), then \(\cocomp{J}\) is an injective object in \(\Dif{A}\).
        \item[2.] \((I,d_I)\) is an injective object in \(\Dif{A}\) if and only if \(\exp(I,d_I)\) is an injective object in \(\Ch{A}\).
    \end{itemize}
\begin{proof}
    \textit{Part 1.} The statement follows from the adjunction pair \((\exp, \cocomp)\), Proposition~\ref{prop: (co)comp are exact} and Lemma~\ref{lem: (F,R) adjunction preservation} part~1.
    
    \textit{Part 2.} By \cite{cartan-homological-1956}*{Prop.\ 2.5 and Prop.\ 2.3}, \((I,d_I)\) is an injective object in \(\Dif{A}\) if and only if \((I,d_I)\) is contractible and \(I\) is an injective object in \(\Mod{A}\). By Proposition~\ref{prop: contractible equivalence} and the description of \(\exp\), this holds if and only if \(\exp(I,d_I)\) is contractible and for each integer \(i\), the evaluation \(\E_i{\exp(I,d_I)}\) is an injective object in \(\Mod{A}\). By \cite{christensen-derived-2024}*{Ex.\ 5.3.4}, this holds if and only if \(\exp(I,d_I)\) is an injective object in \(\Ch{A}\).
\end{proof}
\end{prop}

\begin{prop}\label{prop: preserves semi-injective}
    Let \((I,d_I)\) be an object in \(\Dif{A}\) and \(J\) be an object in \(\Ch{A}\). The following statements hold:
    \begin{itemize}
        \item[1.] If \(J\) is a semi-injective object in \(\Ch{A}\), then \(\cocomp{J}\) is a semi-injective object in \(\Dif{A}\).
        \item[2.] If \((I,d_I)\) is a semi-injective object in \(\Dif{A}\), then \(\exp(I,d_I)\) is a semi-injective object in \(\Ch{A}\).
    \end{itemize}
\begin{proof}
    \textit{Part 1.} The statement follows from the adjunction pair \((\exp, \cocomp)\), Proposition~\ref{prop: (co)comp are exact}, Corollary~\ref{cor: preserves acyclicity and quasi-iso} part~1 and Lemma~\ref{lem: (F,R) adjunction preservation} part~4.
    
    \textit{Part 2.} The statement follows from the adjunction pair \((\comp, \exp)\), Proposition~\ref{prop: (co)comp are exact}, Corollary~\ref{cor: preserves acyclicity and quasi-iso} part~2 and Lemma~\ref{lem: (F,R) adjunction preservation} part~4.
\end{proof}
\end{prop}

\begin{rem}\label{rem: preservation of proj and semi-proj also hold}
    The dual statements of Proposition~\ref{prop: preserves injective} and Proposition~\ref{prop: preserves semi-injective} also hold by their dual arguments. However, we will not use them here.
\end{rem}

\section{The differential tensor and Hom functors}\label{sec: diff tensor and hom}

    Fix an object \(M\) in \(\Mod({A^\op \otimes  B})\).
    
    \paragraph{\bf Differential tensor product.} For an object \((X,d_X)\) in \(\Dif{A}\), we define the \textit{differential tensor product} of \(M\) and \((X,d_X)\) as the object
    \[
    (X,d_X)\boxtimes_A M = (X\otimes_A M,\: d_X \otimes_A M),
    \]
    in \(\Dif{B}\). This assignment augments to an additive functor \(\Dif{A} \xrightarrow[]{-\boxtimes_A M} \Dif{B}\). This functor is related to the functor defined in \cite{avramov-class-2007}*{Sec.\ 1.9}.
    
    \paragraph{\bf Differential Hom.} For an object \((Y,d_Y)\) in \(\Dif{B}\), we define the \textit{differential Hom} of \(M\) and \((Y,d_Y)\) as the object
    \[
    \Dhom_B(M,(Y,d_Y)) = (\Hom_{\Mod{B}}(M,Y),\: \Hom_{\Mod{B}}(M,d_Y))
    \]
    in \(\Dif{A}\). This assignment augments to an additive functor \(\Dif{B} \xrightarrow[]{\Dhom_B(M,-)} \Dif{A}\). These functors form an adjunction pair:
    
\begin{prop}\label{prop: differential tensor Hom adjunction}
    Let \(M\) be an object in \(\Mod({A^\op \otimes  B})\). There is an adjunction isomorphism
    \[
    \Hom_{\Dif{B}}((X,d_X)\boxtimes_{A} M, (Y,d_Y)) \cong \Hom_{\Dif{A}}((X,d_X),\Dhom_B(M,(Y,d_Y)))
    \]
    of objects in \(\Mod{\bbk}\), where \((X,d_X)\) is an object in \(\Dif{A}\) and \((Y,d_Y)\) an object in \(\Dif{B}\). This adjunction isomorphism is given by the canonical currying operation
    \[
    f \mapsto [x \mapsto f(x \otimes -)].
    \]
\begin{proof}
    Let \(X\otimes_A M \xrightarrow[]{f} Y\) be a morphism in \(\Mod{B}\). Then the following diagram on the left is commutative if and only if the following diagram on the right is commutative: 
    \[
    \begin{tikzcd}[column sep = 1.5cm]
        X\otimes_A M \arrow[d, "f"'] \arrow[r, "d_X \otimes_A M"] & X\otimes_A M \arrow[d, "f"] \\
        Y \arrow[r, "d_Y"]                                        & Y                          
    \end{tikzcd}
    \begin{tikzcd}[column sep = 2.5cm]
        X \arrow[r, "d_X"] \arrow[d, "{[x \mapsto f(x\otimes -)]}"']              & X \arrow[d, "{[x \mapsto f(x\otimes -)]}"] \\
        {\Hom_{\Mod{B}}(M,Y)} \arrow[r, "{\Hom_{\Mod{B}}(M,d_Y)}"] & {\Hom_{\Mod{B}}(M,Y)}.
    \end{tikzcd}
    \]
    This means that the Tensor-Hom Adjunction
    \[
    \Hom_{\Mod{B}}(X\otimes_{A} M, Y) \cong \Hom_{\Mod{A}}(X,\Hom_{\Mod{B}}(M,Y))
    \]
    restricts to the desired adjunction isomorphism.
\end{proof}
\end{prop}

    The following result allows us to make new semi-injective objects from old semi-injective objects.

\begin{prop}\label{prop: semi-injective differential Hom}
    Let \(M\) be an object in \(\Mod(A^{\op} \otimes  B)\) and \((I,d_I)\) be an object in \(\Dif{B}\). Assume \(M\) is a projective object in \(\Mod(A^{\op})\) and in \(\Mod{B}\). If \((I,d_I)\) is semi-injective, then \(\Dhom_B(M,(I,d_I))\) is semi-injective in \(\Dif{A}\).
\begin{proof}
    The functors \(\Ch{A} \xrightarrow{-\otimes_{A}M} \Ch{B}\) and \(\Ch{B} \xrightarrow{\Hom_B(M,-)} \Ch{A}\) are exact and preserve acyclic objects (see \cite{christensen-derived-2024}*{Prop.\ 5.2.10, Prop.\ 5.4.9 and Cor.\ 5.4.10}). Consequently, the functors \(-\boxtimes_{A}M\) and \(\Dhom_B(M,-)\) are exact and preserve acyclic objects. The statement follows by applying Lemma~\ref{lem: (F,R) adjunction preservation} to the adjunction pair \((-\boxtimes_A M, \Dhom_B(M,-))\) seen in Proposition~\ref{prop: differential tensor Hom adjunction}.
\end{proof}
\end{prop}

\begin{rem}
    The dual statement of Proposition~\ref{prop: semi-injective differential Hom} concerning the differential tensor product is also true by the dual argument.
\end{rem}

    For the next result, consider the tensor product functor \(\otimes_A\) and the internal Hom functor \(\Hom_B\) on complexes.

\begin{lem}\label{lem: (co)composing tensor Hom}
    Let \(M\) be an object in \(\Mod({A^\op \otimes  B})\). The following statements hold:
    \begin{itemize}
        \item [1.] For each object \(Z\) in \(\Ch{B}\), there is a natural isomorphism
        \[
        \Dhom_B(M,\cocomp{Z}) \cong \cocomp\CHom_B(\S_0{M},Z).
        \]
        \item [2.] For each object \(Y\) in \(\Ch{A}\), there is a natural isomorphism
        \[
        \comp{Y} \boxtimes_A M \cong \comp(Y\otimes_A \S_0{M}).
        \]
    \end{itemize}
\begin{proof}
    \textit{Part 1:} Let \((W,d_W)\) be an object in \(\Dif{A}\). We have the following composition of natural isomorphisms:
    \begin{align*}
        \Hom_{\Dif{A}}((W,d_W),\cocomp\CHom_B(\S_0M,Z)) &\cong \Hom_{\Ch{A}}(\exp(W,d_W), \CHom_B(\S_0{M},Z))\\
                                                     &\cong \Hom_{\Ch{B}}(\exp(W,d_W)\otimes_A \S_0{M}, Z)\\
                                                     &\cong \Hom_{\Ch{B}}(\exp((W,d_W)\boxtimes_A M), Z)\\
                                                     &\cong \Hom_{\Dif{B}}((W,d_W)\boxtimes_A M, \cocomp Z)\\
                                                     &\cong \Hom_{\Dif{A}}((W,d_W), \Dhom_B(M,\cocomp Z)),
    \end{align*}
    where the first and fourth isomorphisms hold by the adjunction pair \((\exp,\cocomp)\), the second isomorphism holds by the Tensor-Hom adjunction, the third isomorphism holds as \(\exp(W,d_W)\otimes_A \S_0{M}\) and \(\exp((W,d_W)\boxtimes_A M)\) are canonically isomorphic and the last isomorphism holds by Proposition~\ref{prop: differential tensor Hom adjunction}. The result now follows from the Yoneda Lemma.

    \textit{Part 2:} This proof follows a similar argument as in part 1.
\end{proof}
\end{lem}

\begin{prop}\label{prop: Isomorphism of diff Hom and comp Hom}
    Let \(J\) be a minimal semi-injective object in \(\Ch{B}\) and \(K\) be a semisimple object in \(\Mod{B}\). For any decomposition \(\cocomp{J} = (I',d_{I'}) \oplus (I'',d_{I''})\) in \(\Dif{B}\), where \((I',d_{I'})\) is a minimal semi-injective object and \((I'',d_{I''})\) an injective object, the following statements hold:
    \begin{itemize}
        \item The object \(\Dhom_B(K,(I'',d_{I''}))\) is zero.
        \item  There is an isomorphism
                \[
                \Dhom_B(K,(I',d_{I'})) \cong \cocomp\CHom_B(\S_0{K},J)
                 \]
               of objects in \(\Dif{B}\).
    \end{itemize}
    In particular, there is an isomorphism
    \[
    \Hom_{\Mod{B}}(K,I') \cong \prod_{i\in\bbZ}\Hom_{\Mod{B}}(K,\E_i{J})
    \]
    of objects in \(\Mod{B}\).
\begin{proof}
    Suppose we have such a decomposition of \(\cocomp{J}\). Applying \(\Dhom_B(K,-)\) to it, we get
    \[
    \Dhom_B(K,(I',d_{I'})) \oplus \Dhom_B(K,(I'',d_{I''})) \cong \Dhom_B(K,\cocomp{J}) \cong \cocomp\CHom_B(\S_0{K},J),
    \]
    where the first isomorphism holds by the additivity of the differential Hom functor and the second isomorphism holds by Lemma~\ref{lem: (co)composing tensor Hom} part 1. As \(J\) is a minimal semi-injective object and \(K\) is semisimple, \(\cocomp\CHom_B(\S_0{K},J)\) has trivial differential (see \cite{christensen-derived-2024}*{Lem.\ 8.2.16}). In particular, \(\Dhom_B(K,(I'',d_{I''}))\) has trivial differential. We will show \(\Dhom_B(K,(I'',d_{I''}))\) is contractible, which then implies it must be zero.
    
    As \((I'',d_{I''})\) is an injective object in \(\Dif{B}\), it is in particular contractible (see \cite{cartan-homological-1956}*{Prop.\ 2.5 and Prop.\ 2.3}). Then Proposition~\ref{prop: contractible equivalence}, gives that the expansion \(\exp(I'',d_{I''})\) is also contractible and, therefore, so is \(\Hom_B(\S_0{K},\exp(I'',d_{I''}))\) (see \cite{christensen-derived-2024}*{Prop.\ 4.3.29}). By Proposition~\ref{prop: contractible equivalence} and as \(\exp\Dhom_B(K,(I'',d_{I''}))\) is canonically isomorphic to \(\Hom_B(\S_0{K},\exp(I'',d_{I''}))\), the object \(\Dhom_B(K,(I'',d_{I''}))\) is contractible. Therefore, \(\Dhom_B(K,(I'',d_{I''}))\) is zero. This gives us the required isomorphism in \(\Dif{B}\) and the following isomorphism in \(\Mod{B}\) is obtained by applying the evaluation functor \(\E_*\) to this isomorphism.
\end{proof}
\end{prop}

\begin{lem}\label{lem: Minimal semi-injective is a direct summan of semiinjective}
    Let \((X,d_X) \xrightarrow[]{} (I,d_I)\) and \((X,d_X) \xrightarrow{} (\Tilde{I},d_{\Tilde{I}})\) be semi-injective resolutions in \(\Dif{A}\). If \((I,d_I)\) is minimal, then \((\Tilde{I},d_{\Tilde{I}})\) has the form \((\Tilde{I},d_{\Tilde{I}})=(I',d_{I'})\oplus (I'',d_{I''})\), where \((I'',d_{I''})\) is an injective object in \(\Dif{A}\) and \((I,d_I)\) is isomorphic to \((I',d_{I'})\) in \(\Dif{A}\).
\begin{proof}
    A special consequence of \cite{holm-minimal-2025}*{Thm.\ B(ii)} is that \((\Tilde{I},d_{\Tilde{I}})\) has the form \((\Tilde{I},d_{\Tilde{I}})= (I',d_{I'}) \oplus (I'',d_{I''})\), with \((I',d_{I'})\) a minimal semi-injective object in \(\Dif{A}\) and \((I'',d_{I''})\) an injective object in \(\Dif{A}\). As it is injective, the object \((I'',d_{I''})\) is contractible and hence acyclic (see \cite{cartan-homological-1956}*{Prop.\ 2.5 and Prop.\ 2.3}). Therefore, the morphism \((X,d_X) \xrightarrow[]{} (\Tilde{I},d_{\Tilde{I}}) \xrightarrow[]{} (I',d_{I'})\), where \((\Tilde{I},d_{\Tilde{I}}) \xrightarrow[]{} (I',d_{I'})\) is the canonical projection, is a minimal semi-injective resolution in \(\Dif{A}\) and hence \((I,d_I)\) and \((I',d_{I'})\) are isomorphic in \(\Dif{A}\) (see \cite{holm-minimal-2025}*{Thm.\ D}).
\end{proof}
\end{lem}

\section{Bass numbers}\label{sec: bass numbers}

    Let \(R\) be a commutative noetherian ring. Matlis' structure theorem \cite{matlis-injective-1958} (also see \cite{christensen-derived-2024}*{Thm.\ C.23}) for injective modules over \(R\) states that every injective \(R\)-module \(E\) decomposes into a coproduct of uniquely determined indecomposable injectives
    \[
    E \cong \coprod_{\fkp\in\Spec{R}}\left(\coprod_{\mu(E_{\fkp})}E_R(R/\fkp)\right),
    \]
    where \(\Spec{R}\) denotes the prime spectrum of the ring \(R\) and the value \(\mu(E_{\fkp})\) is the dimension of \(\Hom_{\Mod{R_{\fkp}}}(\kappa(\fkp),E_{\fkp})\) when viewed as a vector space over \(\kappa(\fkp)=R_{\fkp}/\fkp R_{\fkp}\).
    
    Further assume that \((R,\fkm,k)\) is a local commutative noetherian ring and let \(E\) be an \(R\)-module. We set
    \[
    \mu(E) = \dim_{k}\Hom_{\Mod{R}}(k,E),
    \]
    so that the quantity \(\mu(E)\) is, by Matlis' structure theorem, the number of direct summand copies of \(E_R(k)\) in \(E\). For a semi-injective object \((I,d_I)\) in \(\Dif{R}\) and for a semi-injective object \(J\) in \(\Ch{R}\), we set
    \[
    \mu_D(I,d_I) \coloneqq \mu(I) \text{\quad and \quad} \mu_C(J) \coloneqq \sum_{i\in\bbZ}\mu(\E_i{J}).
    \]
    This makes sense by \cite{holm-qshaped-2024}*{Prop.\ 3.2(b)} for example. Let \(M\) be an object in \(\mod{R}\) and suppose \(M \xrightarrow[]{}J\) is a minimal semi-injective resolution in \(\Ch{R}\). In this case, for an integer \(i\), the quantity \(\mu(\E_i{J})\) is called the \textit{\(i\)-th Bass number} of \(M\). Since the Hom complex \(\Hom_R(k,J)\) has trivial differentials (see \cite{christensen-derived-2024}*{Lem.\ 8.2.16}), the \(i\)-th Bass number can be expressed as
    \[
    \mu(\E_i{J}) = \dim_{k}\Ext^i_{R}(k,M),
    \]
    which is finite (see \cite{christensen-derived-2024}*{Prop.\ 12.3.32}).

\begin{prop}\label{prop: Equal bass numbers}
   Let \((R,\fkm,k)\) be a local commutative noetherian ring and \(J\) be a minimal semi-injective object in \(\Ch{R}\). For any decomposition \(\cocomp{J}= (I',d_{I'}) \oplus (I'',d_{I''})\) in \(\Dif{R}\), where \((I',d_{I'})\) is a minimal semi-injective object and \((I'',d_{I''})\) an injective object, the following equations hold:
   \begin{itemize}
       \item \(\mu_D(I'',d_{I''})=0\).
       \item \(\mu_D(I',d_{I'})=\mu_C(J)\).
   \end{itemize}
\begin{proof}
    This follows by considering the (semi)simple object \(k\) in \(\Mod(k\otimes_{\bbZ}R)\) and applying Proposition~\ref{prop: Isomorphism of diff Hom and comp Hom}.
\end{proof}
\end{prop}

    The following is a differential module version of a classical commutative algebra result (see \cite{bass-ubiquity-1963}*{Lem.\ 2.7} for a classical reference and \cite{christensen-derived-2024}*{Thm.\ 16.4.8} for a modern reference).

\begin{thm}\label{thm: finite injective dimension IFF bass number finite}
    Let \((R,\fkm,k)\) be a local commutative noetherian ring and let \(M\) be an object in \(\mod{R}\) with a minimal semi-injective resolution \(\S_*M \xrightarrow[]{} (I,d_I)\) in \(\Dif{R}\). Then \(M\) has finite injective dimension if and only if \(\mu_D(I,d_I)\) is finite.
\begin{proof}
    Let \(\S_0{M} \xrightarrow[]{}J\) in \(\Ch{R}\) be a minimal semi-injective resolution in \(\Ch{R}\) and consider the compression \(\S_*{M}\cong\cocomp\S_0{M} \xrightarrow[]{}\cocomp{J}\). By Corollary~\ref{cor: preserves acyclicity and quasi-iso} part 6 and Proposition~\ref{prop: preserves semi-injective} part 1, \(\S_*{M}\cong\cocomp\S_0{M} \xrightarrow[]{}\cocomp{J}\) is a semi-injective resolution in \(\Dif{R}\). Lemma~\ref{lem: Minimal semi-injective is a direct summan of semiinjective} gives us a decomposition \(\cocomp{J} \cong (I,d_{I}) \oplus (I'',d_{I''})\), where \((I'',d_{I''})\) is an injective object in \(\Dif{R}\). Proposition~\ref{prop: Equal bass numbers} gives us that \(\mu_D(I,d_I) = \mu_C(J)\). Hence, \(\mu_D(I,d_I)\) is finite if and only if \(\mu_C(J)\) is finite. The equivalence follows since \(M\) has finite injective dimension if and only if \(\mu_C(J)\) is finite (see \cite{christensen-derived-2024}*{Thm.\ 16.4.8}).
\end{proof}
\end{thm}

     Bass asked whether the existence of a finitely generated module of finite injective dimension suffices to characterise local Cohen–Macaulay rings (see \cite{bass-ubiquity-1963}*{Rmk.\ on pg.\ 14}). It turned out to be true. As a direct consequence of this fact, the homological algebra of differential modules can also characterise local Cohen–Macaulay rings.

\begin{cor}[Bass' Question for differential modules]\label{thm: Bass' question}
    A local commutative noetherian ring \((R,\fkm,k)\) is Cohen-Macaulay if and only if there exists a nonzero object \(M\) in \(\mod{R}\) whose minimal semi-injective resolution \(\S_*M \xrightarrow[]{} (I,d_I)\) in \(\Dif{R}\) has that \(\mu_D(I,d_I)\) is finite.

\begin{proof}
    The equivalence follows from \cite{christensen-derived-2024}*{Ex.\ 18.5.9}, \cite{christensen-derived-2024}*{Cor.\ 18.5.8} and Theorem~\ref{thm: finite injective dimension IFF bass number finite}.
\end{proof}
\end{cor}

\section{Preservation of minimal semi-injectives}\label{sec: preservation of min semi-inj}

    It is perhaps unsurprising that the expansion functor interacts well with the notion of minimality:

\begin{prop}\label{prop: preserves minimal semi-injective}
    Let \((I,d_I)\) be a semi-injective object in \(\Dif{A}\). Then \((I,d_I)\) is a minimal semi-injective object in \(\Dif{A}\) if and only if \(\exp(I,d_I)\) is a minimal semi-injective object in \(\Ch{A}\).
\begin{proof}
    \((\xRightarrow[]{\text{only if}}):\) By Proposition~\ref{prop: preserves semi-injective} part 2, the expansion \(\exp(I,d_I)\) is semi-injective object in \(\Ch{A}\). We use the characterisation found in \cite{holm-minimal-2025}*{Thm.\ 3.1(vi)}. Let \(i\) be an integer, \(M\) an object in \(\Mod{A}\) and assume that the morphism
    \[
    \begin{tikzcd}
            \F_i M \arrow[d, "\hat f\: ="] & \cdots \arrow[r] & 0 \arrow[r] \arrow[d] & M \arrow[r, "\id_M"] \arrow[d, "f_0"] & M \arrow[r] \arrow[d, "f_1"] & 0 \arrow[r] & \cdots \\
            \exp(I,d_I)                    & \cdots \arrow[r] & I \arrow[r, "d_I"]      & I \arrow[r, "d_I"]           & I \arrow[r, "d_I"]            & I \arrow[r] & \cdots,
    \end{tikzcd}
    \]
    is a monomorphism in \(\Ch{A}\). As \(\hat{f}\) intertwines the differentials, we must have \(f_1=d_{I}f_0\). We are done if we show that \(M\) is zero. Consider the morphism
    \[
    \begin{tikzcd}[row sep = 2.5em, column sep = 3em]
            \F_* M \arrow[d, "f^\dagger\: ="] && M\oplus M \arrow[r, "{\begin{pmatrix}0 \amsamp 0 \\ \id_M \amsamp0 \end{pmatrix}}"] \arrow[d, "{\begin{pmatrix}f_0 \amsamp d_{I}f_0 \end{pmatrix}}"'] & M\oplus M \arrow[d, "{\begin{pmatrix}f_0 \amsamp d_{I}f_0 \end{pmatrix}}"] \\
            (I,d_I)                                    && I \arrow[r, "d_I"]                                                                                                  & I                                                      
    \end{tikzcd}
    \]
    in \(\Dif{A}\). As \((I,d_I)\) is a minimal semi-injective object in \(\Dif{A}\), it suffices to show that \(f^\dagger\) is a monomorphism (see \cite{holm-minimal-2025}*{Thm.\ 3.1(vi)}). To this end, let \(x=(m,n)\) be an element in \(M\oplus M\) and suppose \(f^\dagger(x)=0\). Applying the differential \(d_I\), the equation \(d_{I}f^\dagger(x)=0\) simplifies to \(d_{I}f_0(m)=0\). As \(\hat{f}\) was assumed to be a monomorphism, \(d_{I}f_0\) is also a monomorphism and we have that \(m=0\). But, then the equation \(f^\dagger(x)=0\) now simplifies to \(d_{I}f_0(n)=0\) and again, we get \(n=0\). Hence \(x=(0,0)\) and \(f^\dagger\) is a monomorphism. Therefore, \(M\) is zero.

    \((\xLeftarrow[]{\text{if}}):\) We use the characterisation found in \cite{holm-minimal-2025}*{Thm.\ 3.1(ii)}. Suppose \((E,d_E)\) is an acyclic subobject of \((I,d_I)\). By Corollary~\ref{cor: preserves acyclicity and quasi-iso} part 1 and as \(\exp\) is an exact functor (see Proposition~\ref{prop: (co)comp are exact} part 1), \(\exp(E,d_E)\) is an acyclic subobject of \(\exp(I,d_I)\). By minimality, \(\exp(E,d_E)\) is zero (see \cite{holm-minimal-2025}*{Thm.\ 3.1(ii)}). As \(\exp\) is faithful, \((E,d_E)\) is also zero.
\end{proof}
\end{prop}

    As the definition of a minimal semi-injective object involves homological and representation theoretic conditions, the question of whether the cocompression functor \(\cocomp{}\) preserves these objects is subtle. Later, we will see in Example~\ref{non-ex: cocomp not preserving minimality} that the minimality of semi-injective objects is not always preserved. Nevertheless, in this section we present some positive results that identify when it does.

\begin{prop}\label{prop: cocompressing minimal semi-injectives of length 2}
    Let \(J\) be an object in \(\Ch{A}\) of the form
    \[
    J = \cdots \xrightarrow[]{} 0 \xrightarrow[]{} J^0 \xrightarrow[]{d^0_J} J^1 \xrightarrow[]{} 0 \xrightarrow[]{} \cdots.
    \]
    If \(J\) is a minimal semi-injective object, then \(\cocomp{J}\) is a minimal semi-injective object in \(\Dif{A}\).
\begin{proof}
    We verify the characterisation of a minimal semi-injective object in \(\Dif{A}\) found in \cite{holm-minimal-2025}*{Thm.\ 3.1(vi)}. By Proposition~\ref{prop: preserves semi-injective} part 1, the cocompression \(\cocomp{J}\) is a semi-injective object in \(\Dif{A}\), and takes the form
    \[
    \cocomp{J}= \left(J^0\oplus J^1 ,{\begin{pmatrix}0 \amsamp 0 \\ d^0_J \amsamp0 \end{pmatrix}}\right).
    \]
    Let \(M\) be an object in \(\Mod{A}\) and assume that the morphism \(f^\dagger\):
    \[
    \begin{tikzcd}[row sep = 4em, column sep = 3em]
        \F_* M \arrow[d, "f^\dagger \:="] && M\oplus M \arrow[r, "{\begin{pmatrix}0 \amsamp 0 \\ \id_M \amsamp0 \end{pmatrix}}"] \arrow[d, "{\begin{pmatrix}f_{00} \amsamp f_{01} \\ f_{10} \amsamp f_{11} \end{pmatrix}}"'] & M\oplus M \arrow[d, "{\begin{pmatrix}f_{00} \amsamp f_{01} \\ f_{10} \amsamp f_{11} \end{pmatrix}}"] \\
            \cocomp{J}                                    && J^0\oplus J^1 \arrow[r, "{\begin{pmatrix}0 \amsamp 0 \\ d^0_J \amsamp0 \end{pmatrix}}"]                                                                                                  & J^0\oplus J^1,
    \end{tikzcd}
    \]
    is a monomorphism in \(\Dif{A}\). As \(f^\dagger\) intertwines with the differentials, it must be a monomorphism of the form
    \[
    M\oplus M \xrightarrow[]{{\begin{pmatrix}f_{00} \amsamp 0 \\ f_{10} \amsamp d^0_Jf_{00} \end{pmatrix}}} J^0 \oplus J^1.
    \]
    In particular, \(d^0_Jf_{00} \colon M \xrightarrow[]{} J^1\) is a monomorphism. But then
    \[
    \begin{tikzcd}
        \cdots \arrow[r] & 0 \arrow[r] \arrow[d] & M \arrow[r, "\id_M"] \arrow[d, "f_{00}"] & M \arrow[r] \arrow[d, "d^0_Jf_{00}"] & 0 \arrow[r] & \cdots    \\
        \cdots \arrow[r] & 0 \arrow[r]           & J^0 \arrow[r, "d^0_J"]         & J^1 \arrow[r]                   & 0 \arrow[r] & {\cdots,}
    \end{tikzcd}
    \]
    is a monomorphism in \(\Ch{A}\) and hence, by minimality of \(J\), the object \(M\) is zero.
\end{proof}
\end{prop}

    Let \(M\) be an object in \(\Mod{A}\) and consider a submodule \(L\) of \(M\). We say \(L\) is \textit{large} (also known as \textit{essential}) if the zero object is the only submodule of \(M\) with a trivial intersection with \(L\). Let \(\Soc M\) denote the socle of \(M\). By \cite{anderson-rings-1992}*{Prop.\ 9.7}, we know that \(\Soc{M}\) is equal to the intersection of all of the large submodules of \(M\).

\begin{prop}\label{prop: cocompressing minimal semi-injectives of fin gen}
    Let \(A\) be a right Artinian \(\bbk\)-algebra and \(J\) be an object in \(\Ch{A}\) of the form
    \[
    J = \cdots \xrightarrow[]{} J^{0} \xrightarrow[]{d^{0}_J} J^1 \xrightarrow[]{d^1_J} J^2 \xrightarrow[]{} \cdots.
    \]
    Suppose \(J\) is bounded. If \(J\) is a minimal semi-injective object, then \(\cocomp{J}\) is a minimal semi-injective object in \(\Dif{A}\).
\begin{proof}
    We verify the characterisation of a minimal semi-injective object in \(\Dif{A}\) found in \cite{holm-minimal-2025}*{Thm.\ 3.1(vi)}. By Proposition~\ref{prop: preserves semi-injective} part 1, the cocompression \(\cocomp{J}\) is a semi-injective object in \(\Dif{A}\). As \(J\) is bounded, we can assume its cocompression is of the form
    \[
       \cocomp{J}= \left( \bigoplus_{i=0}^{n} J^i, d_J \coloneqq \bigoplus_{i=0}^{n}d^i_J\right),
    \]
    for \(n\) some natural number. Let \(M\) be an object in \(\Mod{A}\) and assume that the morphism \(f^\dagger\):
    \[
    \begin{tikzcd}[row sep = 5ex]
        F_* M \arrow[d, "f^\dagger \: ="] & M\oplus M \arrow[r, "{\begin{pmatrix}0 \amsamp 0 \\ \id_M \amsamp0 \end{pmatrix}}"] \arrow[d, "{\begin{pmatrix}f_0 \amsamp f_1 \end{pmatrix}}"'] & M\oplus M \arrow[d, "{\begin{pmatrix}f_0 \amsamp f_1 \end{pmatrix}}"] \\
            \cocomp{J}                                    & \bigoplus_{i=0}^{n} J^i \arrow[r, "d_J"]                                                                                                  & \bigoplus_{i=0}^{n} J^i                                                   
    \end{tikzcd}
    \]
    is a monomorphism in \(\Dif{A}\). As \(f^\dagger\) intertwines with the differentials, it must be of the form
    \[
    M\oplus M \xrightarrow[]{{\begin{pmatrix}f_0 \amsamp d_Jf_0 \end{pmatrix}}} \bigoplus_{i=0}^{n} J^i.
    \]
    In particular, the component \(d_Jf_0\) is also a monomorphism in \(\Mod{A}\). The socle of an object in \(\Mod{A}\) is equal to the intersection of all of its large submodules (see \cite{anderson-rings-1992}*{Prop.\ 9.7}) and as \(J\) is minimal, the submodule \(\Ker{d^i_J}\) is large in \(J_i\) (see \cite{christensen-derived-2024}*{Theorem~B.21}) and therefore, the direct sum \(\Ker{d_J} \cong \bigoplus_{i=0}^n\Ker{d^i_J}\) is large in \(\bigoplus_{i=0}^{n}J_i\) (see \cite{goodearl-introduction-1989}*{Prop.\ 5.6(d)}). This implies that \(\Soc(\bigoplus_{i=0}^{n}J_i)\) is contained in \(\Ker d_J\). Now, the component \(f_0\) maps \(\Soc{M}\) into \(\Soc(\bigoplus_{i=0}^{n}J_i)\) and therefore, the composition \(d_Jf_0\) maps \(\Soc{M}\) to zero. As \(d_Jf_0\) is a monomorphism, \(\Soc{M}\) is zero. But the only object in \(\Mod{A}\) with a trivial socle is the zero object (see \cite{auslander-representation-1997}*{Prop.\ 4.1(a)}). Note that the proof of \cite{auslander-representation-1997}*{Prop.\ 4.1(a)} remains true when the finitely generated assumption and the assumption that the ring is an algebra over a field is dropped (see \cite{auslander-representation-1997}*{Prop.\ 3.1(a)}). Hence, \(M\) is zero.
\end{proof}
\end{prop}

\begin{rem}
    Recall Proposition~\ref{prop: Equal bass numbers}. It tells us that over a local commutative noetherian ring, all the indecomposable injective summands \(E(k)\) appearing in a minimal semi-injective complex \(J\), when cocompressed by \(\cocomp{}\), must land in the minimal semi-injective direct summand \((I',d_{I'})\) of the decomposition \(\cocomp{J}=(I',d_{I'}) \oplus (I'',d_{I''})\). We will now use this fact to show another instance for when the cocompression functor \(\cocomp\) preserves minimal semi-injective objects.
\end{rem}

The next preservation result applies to many rings, for example:
\begin{itemize}
        \item \(\bbZ/p^n\bbZ\) for \(p\) a prime.
        \item \(k[x]/(f^n)\) for \(f\) an irreducible polynomial.
        \item \(k[x^2,x^3]/(x^3)\).
\end{itemize}

\begin{prop}\label{prop: cocomp preserves local noetherian krull dimension one}
    Let \((R,\fkm,k)\) be a local commutative noetherian ring of Krull dimension zero (or equivalently, \(R\) is a local commutative Artinian ring) and let \(J\) be an object in \(\Ch{R}\). If \(J\) is a minimal semi-injective object, then \(\cocomp{J}\) is a minimal semi-injective object in \(\Dif{R}\).
\begin{proof}
    By Proposition~\ref{prop: preserves semi-injective} part 1, the cocompression \(\cocomp{J}\) is a semi-injective object in \(\Dif{R}\). By \cite{holm-minimal-2025}*{Thm.\ B(ii)}, it is of the form \(\cocomp{J}= (I',d_{I'}) \oplus (I'',d_{I''})\), where \((I',d_{I'})\) is minimal semi-injective and \((I'',d_{I''})\) injective. It suffices to show that \((I'',d_{I''})\) is zero. As \(R\) is of Krull dimension zero, \(E(k)\) is the unique indecomposable injective object in \(\Mod{R}\). By Matlis' structure theorem \cite{matlis-injective-1958} (also see \cite{christensen-derived-2024}*{Thm.\ C.23}), each injective object \(E\) in \(\Mod{R}\) is a coproduct of \(\mu(E)\) many summands isomorphic to \(E(k)\). By Proposition~\ref{prop: Equal bass numbers}, the equation \(\mu_D(I'',d_{I''})=0\) holds and thus, \((I'',d_{I''})\) must be zero.
\end{proof}
\end{prop}

\section{Examples}\label{sec: examples}
    
\begin{ex}[Hereditary rings]\label{ex: hereditary rings}
    Considering Proposition~\ref{prop: cocompressing minimal semi-injectives of length 2}, we take \(A\) to be the ring of integers \(\bbZ\) and let \(I\) be the minimal semi-injective resolution of \(\S_0\bbZ\) in \(\Ch{\bbZ}\):
    \[
    I = \cdots \xrightarrow[]{} 0 \xrightarrow[]{} \bbQ \xrightarrow[]{\pi} \bbQ/\bbZ \xrightarrow[]{} 0 \xrightarrow[]{} \cdots.
    \]
    Its cocompression
    \[
    \cocomp{I}= \left(\bbQ \oplus \bbQ/\bbZ, \begin{pmatrix}0 \amsamp 0 \\ \pi \amsamp0 \end{pmatrix} \right),
    \]
    is therefore a minimal semi-injective resolution of \(\S_*\bbZ\) in \(\Dif{\bbZ}\). By Proposition~\ref{prop: preserves minimal semi-injective}, the expansion
    \[
    \exp\cocomp{I}= \cdots \xrightarrow[]{} \bbQ \oplus \bbQ/\bbZ \xrightarrow[]{{\begin{pmatrix}0 \amsamp 0 \\ \pi \amsamp0 \end{pmatrix}}} \bbQ \oplus \bbQ/\bbZ \xrightarrow[]{{\begin{pmatrix}0 \amsamp 0 \\ \pi \amsamp0 \end{pmatrix}}} \bbQ \oplus \bbQ/\bbZ \xrightarrow[]{{\begin{pmatrix}0 \amsamp 0 \\ \pi \amsamp0 \end{pmatrix}}} \bbQ \oplus \bbQ/\bbZ \xrightarrow[]{} \cdots
    \]
    is a minimal semi-injective object in \(\Ch{A}\).
\end{ex}

\begin{ex}[Infinite global dimension]\label{ex: Finite dimensional algebras 2}
    Considering Proposition~\ref{prop: cocomp preserves local noetherian krull dimension one}, we take \((R,\fkm)\) to be \((\bbZ/4\bbZ, 2\bbZ/4\bbZ)\) or the dual numbers \((k[x]/(x^2), (x)/(x^2))\) over a field \(k\). Let \(I\) be the minimal semi-injective resolution of the residue field \(\S_0(R/\fkm)\) in \(\Ch{R}\):
    \[
    I = \cdots \xrightarrow[]{} 0 \xrightarrow[]{} R \xrightarrow[]{\delta} R \xrightarrow[]{\delta} R \xrightarrow[]{} \cdots.
    \]
    Its cocompression is of the form
    \[
    \cocomp{I} = \left(\prod_{i\in\bbN} R, \:\overline{\delta} \right),
    \]
     where the differential \(\overline{\delta}\) is given by \((t_0,t_1,t_2,\ldots) \mapsto (0,\delta t_0,\delta t_1,\ldots)\), and it is a minimal semi-injective resolution of \(\S_*(R/\fkm)\) in \(\Dif{R}\). Using Proposition~\ref{prop: preserves minimal semi-injective}, the expansion
    \[
    \exp\cocomp{I}= \cdots \xrightarrow[]{}\prod_{i\in\bbN} R \xrightarrow[]{\overline{\delta}}\prod_{i\in\bbN} R \xrightarrow[]{\overline{\delta}}\prod_{i\in\bbN} R \xrightarrow[]{} \cdots
    \]
    is a minimal semi-injective object in \(\Ch{R}\).
\end{ex}

    We will now see that the cocompression functor does not always preserve minimal semi-injective objects.

\begin{ex}[Non-example hereditary rings]\label{non-ex: cocomp not preserving minimality}
    Let \(A\) be the ring of integers \(\bbZ\) and let \(J\) be the minimal semi-injective resolution of \(\S_0\bbZ\) as an object in \(\Ch{\bbZ}\):
    \[
    J = \cdots \xrightarrow[]{} 0 \xrightarrow[]{} \bbQ \xrightarrow[]{\pi} \bbQ/\bbZ \xrightarrow[]{} 0 \xrightarrow[]{} \cdots.
    \]
    We consider the following object in \(\Ch{A}\):
    \begin{equation}\label{eqn: 1}
      \prod_{i\geq 0}\Sigma^{-i}(J) = \cdots \xrightarrow[]{} 0 \xrightarrow[]{} \bbQ \xrightarrow[]{{\begin{pmatrix}0 \\ \pi \end{pmatrix}}} \bbQ \oplus \bbQ/\bbZ \xrightarrow[]{{\begin{pmatrix}0 \amsamp 0 \\ \pi \amsamp0 \end{pmatrix}}} \bbQ \oplus \bbQ/\bbZ \xrightarrow[]{{\begin{pmatrix}0 \amsamp 0 \\ \pi \amsamp0 \end{pmatrix}}} \bbQ \oplus \bbQ/\bbZ \xrightarrow[]{} \cdots.
    \end{equation}
    It is a semi-injective object in \(\Ch{\bbZ}\) as it is bounded on the left and consists of injective objects in \(\Mod{A}\) (see \cite{christensen-derived-2024}*{Ex.\ 5.3.12}, but note that homological grading is used). It is minimal as the kernel of each differential is large (see \cite{holm-minimal-2025}*{Thm.\ 3.1(vii)}). By Proposition~\ref{prop: preserves semi-injective} part 1, the cocompression
    \begin{equation}\label{eqn: 3}
        \cocomp\left(\prod_{i\geq 0}\Sigma^{-i}(J)\right)
    \end{equation}
    is a semi-injective object in \(\Dif{\bbZ}\). We will show that it is not minimal. By Proposition~\ref{prop: preserves minimal semi-injective}, it is minimal if and only if its expansion
    \begin{equation}\label{eqn: 2}
        \exp\cocomp\left(\prod_{i\geq 0}\Sigma^{-i}(J)\right)
    \end{equation}
    is minimal. As \eqref{eqn: 2} is periodic, minimality holds if and only if the kernel of the differential of \eqref{eqn: 3} is large (see \cite{holm-minimal-2025}*{Thm.\ 3.1(vii)}). The kernel of the differential of \eqref{eqn: 3} is isomorphic to the product of the kernels of the differentials in \eqref{eqn: 1}. So, we need to show that the submodule 
    \[
    \Ker \cocomp\left(\prod_{i\geq 0}\Sigma^{-i}(J)\right) \cong \bbZ \oplus \prod_{i\geq1} (\bbZ \oplus \bbQ/\bbZ) \subseteq \bbQ \oplus \prod_{i\geq1} (\bbQ \oplus \bbQ/\bbZ) \cong \cocomp\left(\prod_{i\geq 0}\Sigma^{-i}(J)\right)
    \]
    is not large. To this end, let \(M\) be the submodule of \(\bbQ \oplus \prod_{i\geq1} (\bbQ \oplus \bbQ/\bbZ)\) generated by the element \(x=(0, (\frac{1}{i},0)_{i\geq1})\). As there is no integer \(n\) such that \(\frac{n}{i}\) is an integer for all \(i\geq 1\), there cannot exist a nonzero element of \(M\) that also lies in \(\bbZ \oplus \prod_{i\geq1} (\bbZ \oplus \bbQ/\bbZ)\). Hence, the kernel is not large.
\end{ex}

\paragraph{\bf Acknowledgement.} I am grateful to my supervisor, Peter J{\o}rgensen, for introducing me to the theory of the \(Q\)-shaped derived category and for his invaluable guidance throughout the development of this paper. I also thank Cyril Matou{\v s}ek for our enjoyable mathematical discussions. I thank the anonymous referee for their careful reading of the manuscript and for their corrections and suggestions. This work was supported by the Independent Research Fund Denmark (grant no. 1026-00050B).

\footnotesize
\bibliographystyle{alpha}
\bibliography{references}

@book{christensen-derived-2024,
	title        = {Derived {C}ategory {M}ethods in {C}ommutative {A}lgebra},
	author       = {Christensen, Lars W. and Foxby, Hans-Bj{\o}rn and Holm, Henrik},
	year         = {2024},
	publisher    = {Springer Cham},
	series       = {Springer Monographs in Mathematics},
	pages        = {xxiii+1119},
	doi          = {10.1007/978-3-031-77453-9},
	edition      = {1st ed.}
}

@article {holm-minimal-2025,
    AUTHOR = {Holm, Henrik and J{\o}rgensen, Peter},
     TITLE = {Minimal semiinjective resolutions in the {$Q$}-shaped derived
              category},
   JOURNAL = {Rev. Mat. Iberoam.},
  FJOURNAL = {Revista Matem\'atica Iberoamericana},
    VOLUME = {41},
      YEAR = {2025},
    NUMBER = {6},
     PAGES = {2309--2334},
       DOI = {10.4171/rmi/1579},
}

@article{holm-qshaped-2024,
	title        = {The {$Q$}-shaped derived category of a ring---compact and perfect objects},
	author       = {Holm, Henrik and J{\o}rgensen, Peter},
	year         = {2024},
	journal      = {Trans. Amer. Math. Soc.},
	volume       = {377},
	number       = {5},
	pages        = {3095--3128},
	doi          = {10.1090/tran/8979},
	fjournal     = {Transactions of the American Mathematical Society}
}

@article{holm-qshaped-2022,
	title        = {The {$Q$}-shaped derived category of a ring},
	author       = {Holm, Henrik and J{\o}rgensen, Peter},
	year         = {2022},
	journal      = {J. Lond. Math. Soc. (2)},
	volume       = {106},
	number       = {4},
	pages        = {3263--3316},
	doi          = {10.1112/jlms.12662},
	fjournal     = {Journal of the London Mathematical Society. Second Series}
}

@article{holm-cotorsion-2019,
	title        = {Cotorsion pairs in categories of quiver representations},
	author       = {Holm, Henrik and J{\o}rgensen, Peter},
	year         = {2019},
	journal      = {Kyoto J. Math.},
	volume       = {59},
	number       = {3},
	pages        = {575--606},
	doi          = {10.1215/21562261-2018-0018},
	fjournal     = {Kyoto Journal of Mathematics}
}

@article{ringel-representations-2017,
	title        = {Representations of quivers over the algebra of dual numbers},
	author       = {Ringel, Claus M. and Zhang, Pu},
	year         = {2017},
	journal      = {J. Algebra},
	volume       = {475},
	pages        = {327--360},
	doi          = {10.1016/j.jalgebra.2016.12.001},
	fjournal     = {Journal of Algebra}
}

@article{avramov-class-2007,
	title        = {Class and rank of differential modules},
	author       = {Avramov, Luchezar L. and Buchweitz, Ragnar-Olaf and Iyengar, Srikanth},
	year         = {2007},
	journal      = {Invent. Math.},
	volume       = {169},
	number       = {1},
	pages        = {1--35},
	doi          = {10.1007/s00222-007-0041-6},
	fjournal     = {Inventiones Mathematicae}
}

@book{maclane-categories-1998,
	title        = {Categories for the working mathematician},
	author       = {Mac Lane, Saunders},
	year         = {1998},
	publisher    = {Springer-Verlag, New York},
	series       = {Graduate Texts in Mathematics},
	volume       = {5},
	pages        = {xii+314},
	edition      = {2nd ed.}
}

@book{anderson-rings-1992,
	title        = {Rings and categories of modules},
	author       = {Anderson, Frank W. and Fuller, Kent R.},
	year         = {1992},
	publisher    = {Springer-Verlag, New York},
	series       = {Graduate Texts in Mathematics},
	volume       = {13},
	pages        = {x+376},
	doi          = {10.1007/978-1-4612-4418-9},
	edition      = {2nd ed.}
}

@book{stenstrom-grundlehren-1975,
	title        = {Rings of quotients},
	author       = {Stenstr{\"o}m, Bo},
	year         = {1975},
	publisher    = {Springer Berlin, Heidelberg},
	series       = {Die Grundlehren der mathematischen Wissenschaften},
	pages        = {viii+309},
	doi          = {10.1007/978-3-642-66066-5},
	subtitle     = {An {I}ntroduction to {M}ethods of {R}ing {T}heory},
	edition      = {1st ed.}
}

@article{matlis-injective-1958,
	title        = {Injective modules over {N}oetherian rings},
	author       = {Matlis, Eben},
	year         = {1958},
	journal      = {Pacific J. Math.},
	volume       = {8},
	pages        = {511--528},
	url          = {http://projecteuclid.org/euclid.pjm/1103039896},
	fjournal     = {Pacific Journal of Mathematics}
}

@article{grothendieck-quelques-1957,
	title        = {Sur quelques points d'alg{\`e}bre homologique},
	author       = {Grothendieck, Alexander},
	year         = {1957},
	journal      = {Tohoku Math. J. (2)},
	volume       = {9},
	pages        = {119--221},
	doi          = {10.2748/tmj/1178244839},
	fjournal     = {The Tohoku Mathematical Journal. Second Series}
}

@book{cartan-homological-1956,
	title        = {Homological algebra},
	author       = {Cartan, Henri and Eilenberg, Samuel},
	year         = {1956},
	publisher    = {Princeton University Press, Princeton, NJ},
	pages        = {xv+390}
}

@book {goodearl-introduction-1989,
    AUTHOR = {Goodearl, Kenneth R. and Warfield, Jr., Robert B.},
     TITLE = {An introduction to noncommutative {N}oetherian rings},
    SERIES = {London Mathematical Society Student Texts},
    VOLUME = {16},
 PUBLISHER = {Cambridge University Press, Cambridge},
      YEAR = {1989},
     PAGES = {xviii+303},
DOI = {10.1017/CBO9780511841699}
}

@book {auslander-representation-1997,
    AUTHOR = {Auslander, Maurice and Reiten, Idun and Smal{\o}, Sverre O.},
     TITLE = {Representation theory of {A}rtin algebras},
    SERIES = {Cambridge Studies in Advanced Mathematics},
    VOLUME = {36},
      NOTE = {Corrected reprint of the 1995 original},
 PUBLISHER = {Cambridge University Press, Cambridge},
      YEAR = {1997},
     PAGES = {xiv+425},
}

@article {bass-ubiquity-1963,
    AUTHOR = {Bass, Hyman},
     TITLE = {On the ubiquity of {G}orenstein rings},
   JOURNAL = {Math. Z.},
  FJOURNAL = {Mathematische Zeitschrift},
    VOLUME = {82},
      YEAR = {1963},
     PAGES = {8--28},
       DOI = {10.1007/BF01112819},
}

@article {rouquier-representation-2006,
    AUTHOR = {Rouquier, Rapha{\"e}l},
     TITLE = {Representation dimension of exterior algebras},
   JOURNAL = {Invent. Math.},
  FJOURNAL = {Inventiones Mathematicae},
    VOLUME = {165},
      YEAR = {2006},
    NUMBER = {2},
     PAGES = {357--367},
       DOI = {10.1007/s00222-006-0499-7},
}

@article {stai-differential-2017,
    AUTHOR = {Stai, Torkil},
     TITLE = {Differential modules over quadratic monomial algebras},
   JOURNAL = {Algebr. Represent. Theory},
  FJOURNAL = {Algebras and Representation Theory},
    VOLUME = {20},
      YEAR = {2017},
    NUMBER = {5},
     PAGES = {1239--1247},
       DOI = {10.1007/s10468-017-9684-1},
}

@article {brown-minimal-2022,
    AUTHOR = {Brown, Michael K. and Erman, Daniel},
     TITLE = {Minimal free resolutions of differential modules},
   JOURNAL = {Trans. Amer. Math. Soc.},
  FJOURNAL = {Transactions of the American Mathematical Society},
    VOLUME = {375},
      YEAR = {2022},
    NUMBER = {10},
     PAGES = {7509--7528},
       DOI = {10.1090/tran/8754},
}

@article {wei-gorenstein-2015,
    AUTHOR = {Wei, Jiaqun},
     TITLE = {Gorenstein homological theory for differential modules},
   JOURNAL = {Proc. Roy. Soc. Edinburgh Sect. A},
  FJOURNAL = {Proceedings of the Royal Society of Edinburgh. Section A.
              Mathematics},
    VOLUME = {145},
      YEAR = {2015},
    NUMBER = {3},
     PAGES = {639--655},
       DOI = {10.1017/S0308210513000541},
}

@article {peskine-dimension-1973,
    AUTHOR = {Peskine, Christian and Szpiro, Lucien},
     TITLE = {Dimension projective finie et cohomologie locale},
SUBTITLE = {{A}pplications {\`a} la d{\'e}monstration de conjectures de {M}.
              {A}uslander, {H}. {B}ass et {A}. {G}rothendieck},
   JOURNAL = {Inst. Hautes {\'E}tudes Sci. Publ. Math.},
  FJOURNAL = {Institut des Hautes {\'E}tudes Scientifiques. Publications
              Math{\'e}matiques},
    VOLUME = {42},
      YEAR = {1973},
     PAGES = {47--119},
       URL = {http://www.numdam.org/item?id=PMIHES_1973__42__47_0},
}

@article {roberts-theoreme-1987,
    AUTHOR = {Roberts, Paul},
     TITLE = {Le th{\'e}or{\`e}me d'intersection},
   JOURNAL = {C. R. Acad. Sci. Paris S{\'e}r. I Math.},
  FJOURNAL = {Comptes Rendus des S{\'e}ances de l'Acad{\'e}mie des Sciences.
              S{\'}erie I. Math{\'e}matique},
    VOLUME = {304},
      YEAR = {1987},
    NUMBER = {7},
     PAGES = {177--180},
}

@book {roberts-multiplicities-1998,
    AUTHOR = {Roberts, Paul C.},
     TITLE = {Multiplicities and {C}hern classes in local algebra},
    SERIES = {Cambridge Tracts in Mathematics},
    VOLUME = {133},
 PUBLISHER = {Cambridge University Press, Cambridge},
      YEAR = {1998},
     PAGES = {xii+303},
       DOI = {10.1017/CBO9780511529986},
}

@article {auslander-homological-1956,
    AUTHOR = {Auslander, Maurice and Buchsbaum, David A.},
     TITLE = {Homological dimension in {N}oetherian rings},
   JOURNAL = {Proc. Nat. Acad. Sci. U.S.A.},
  FJOURNAL = {Proceedings of the National Academy of Sciences of the United
              States of America},
    VOLUME = {42},
      YEAR = {1956},
     PAGES = {36--38},
       DOI = {10.1073/pnas.42.1.36},
}

@inproceedings {serre-dimension-1956,
    AUTHOR = {Serre, Jean-Pierre},
     TITLE = {Sur la dimension homologique des anneaux et des modules
              noeth\'eriens},
 BOOKTITLE = {Proceedings of the international symposium on algebraic number
              theory, {T}okyo \& {N}ikko, 1955},
     PAGES = {175--189},
 PUBLISHER = {Science Council of Japan, Tokyo},
      YEAR = {1956},
}

\textsc{Department of Mathematics, Aarhus University, Ny Munkegade 118, 8000 Aarhus C, Denmark}.\\
Email address: \href{mailto:david.nkansah@math.au.dk}{david.nkansah@math.au.dk}

\end{document}